\documentclass[12pt]{article}
\usepackage{amsfonts,amsmath,amssymb,amscd}
\newtheorem{definition}{Definition}[section]

\newtheorem{proposition}[definition]{Proposition}
\newtheorem{theorem}[definition]{Theorem}

                    \def\cF{{\cal F}}
                    \def\cI{{\cal I}}

\def\cS{{\cal S}}                    \def\cU{{\cal U}}

\def\fB{{\mathfrak B}}
\def\fC{{\mathfrak C}}

\def\fa{{\mathfrak a}}
\def\fg{{\mathfrak g}}

\newcommand{\CC}{{\mathbb C}}
\newcommand{\II}{{\mathbb I}}

\newcommand{\ZZ}{{\mathbb Z}}

\newcommand{\un}{\mbox{1\hspace{-1mm}I}}

\newcommand{\finproof}{{\hfill \rule{5pt}{5pt}}}

\def\qmbox#1{\qquad\mbox{#1}\quad}

\def\QTHA{Quasitriangular Hopf Algebra (QTHA)\def\QTHA{QTHA}}
\def\QTQHA{Quasitriangular Quasi-Hopf Algebra (QTQHA)\def\QTQHA{QTQHA}}

\begin{document}
\pagestyle{empty}


\begin{center}

{\Large \textsf{$R$-matrix presentation for (super) Yangians
    $Y(\fg)$}} 

\vspace{10mm}

{\large D. Arnaudon$^a$, J. Avan$^b$, N.~Cramp\'e$^a$,
  L. Frappat$^{ac}$, {E}. Ragoucy$^a$}

\vspace{10mm}

\emph{$^a$ Laboratoire d'Annecy-le-Vieux de Physique Th{\'e}orique}

\emph{LAPTH, CNRS, UMR 5108, Universit{\'e} de Savoie}

\emph{B.P. 110, F-74941 Annecy-le-Vieux Cedex, France}

\vspace{7mm}

\emph{$^b$ Laboratoire de Physique Th{\'e}orique et Hautes {\'E}nergies}

\emph{LPTHE, CNRS, UMR 7589, Universit{\'e}s Paris VI/VII}

\emph{4, place Jussieu, B.P. 126, F-75252 Paris Cedex 05, France}

\vspace{7mm}

\emph{$^c$ Member of Institut Universitaire de France}

\end{center}

\vfill
\vfill

\begin{abstract}
We give a unified RTT presentation of (super)-Yangians $Y(\fg)$ for
$\fg=so(n)$, $sp(2n)$ and $osp(m|2n)$. 
\end{abstract}

\vfill
MSC number: 81R50, 17B37
\vfill

\rightline{LAPTH-878/01}
\rightline{PAR-LPTHE 01-70}
\rightline{math.QA/0111325}

\baselineskip=16pt

\newpage

\pagestyle{plain}
\setcounter{page}{1}

\section{Introduction}
\setcounter{equation}{0}

The Yangian $Y(\fa)$ based on a simple Lie algebra $\fa$ is defined
\cite{Dri85,Dri86} as the
homogeneous quantisation of the algebra $\fa[u] = \fa \otimes
\CC[u]$ endowed with
its standard bialgebra structure, where $\CC[u]$ is the ring of
polynomials in the 
indeterminate $u$. There exists for the Yangian $Y(\fa)$ three
different realisations, due to Drinfel'd \cite{Dri85,Dri86,Dri88}. 
In the first realisation the Yangian is
generated by the elements $J^a_0$ of the Lie algebra and a second set
of  generators $J^a_1$ in one-to-one correspondence with $J^a_0$
realising a representation space thereof. The second
realisation is given in terms of generators and relations similar 
to the description of a loop algebra as a space of maps. 
However in this realisation no explicit formula for the 
comultiplication is known in general, except in the $sl(2)$ case
\cite{Molev2001}. The third realisation uses the 
Faddeev--Reshetikhin--Takhtajan (FRT) 
formalism \cite{FRT}, but it is only established in
the cases of classical Lie algebras.

The FRT formalism is also used as the original definition of the
super Yangian $Y(gl(M|N))$~\cite{Naz91,Zha95}. 
The purpose of this paper is to define the Yangian for the 
orthosymplectic Lie superalgebras via the FRT formalism. 
As a by-product, we exhibit a unified construction 
which encompasses the three cases $\fg=so(M)$, $\fg=sp(N)$ and
$\fg=osp(M|N)$. A key feature in this procedure is the explicit
expression of a 
``quantum determinant''-like central element which coincides with that
given by Drinfel'd in the $\fg=so(M)$ case \cite{Dri85}.

Note that a first attempt for an FRT 
formulation of Yangians based on $so(M)$ and $sp(N)$ was done by 
Olshanski \textit{et al.}~\cite{olsh,MNO}. However, it led to the notion of 
twisted Yangians, which indeed are deformations of loop algebras on 
$so(M)$ and $sp(N)$, but appear as Hopf coideals rather than Hopf 
algebras. The same feature holds for twisted super Yangians, 
corresponding to $osp(M|N)$ superalgebras~\cite{briot}.

Known rational solutions of the Yang--Baxter equation involve
$R$-matrices of the form \emph{(i)} $R(u) = \II + \frac{P}{u}$ and 
\emph{(ii)} $R(u) = \II + \frac{P}{u} - \frac{K}{u+\kappa}$  
\cite{ABZABZ78,KS,Isa}. 
The first case , where $P$ is defined as the (super)-permutation map, is
known to define the Yangians $Y(sl(N))$ and $Y(sl(M|N))$ via the FRT 
formalism \cite{FRT,KR}. In the case \emph{(ii)}, 
$K$ is a partial (super)-transposition of $P$. Some $R$-matrices of
this type occur
as factorised $S$-matrices of quantum field models in two dimensions
exhibiting the $so(M)$ symmetry \cite{ABZABZ78}.

We will show that the $R$-matrix 
$R(u) = \II + \frac{P}{u} - \frac{K}{u+\kappa}\;$ can be used to define
the Yangian $Y(\fg)$ within the RTT formalism, for 
$\fg=osp(M|N)$ ($N$ even) as well as for the
cases of $\fg=so(M)$ or $\fg=sp(N)$ ($N$ even). 
We prove that the algebra defined this way 
is indeed a quantisation of $\fg[u]$ endowed with its
canonical bialgebra structure.

The letter is organised as follows. In Section \ref{sect:general},
after some definitions, we introduce for each $\fg$ a rational
$R$-matrix expressed 
in terms of the (super)-permutation and of its partial transposition. 
We check that it satisfies the (super) Yang--Baxter equation in all
cases. In Section 
\ref{sect:yangians} we define a (super)-algebra through the RTT
formalism. We establish that the quotient of this algebra by 
the quantum determinant-like central 
element is the Yangian~$Y(\fg)$, as defined
in \cite{Dri85,Dri86}.

\section{General setting\label{sect:general}}
\setcounter{equation}{0}

Let $gl(M|N)$ be the $\ZZ_{2}$-graded algebra of $(M+N) \!\times\!
(M+N)$ matrices $X_{ij}$. Let $\theta_0=\pm 1$. The
$\ZZ_{2}$-gradation is defined by
$(-1)^{[i]} = \theta_0$ if $1 \le i \le M$ and $(-1)^{[i]} = -\theta_0$
if $M+1 \le i \le M+N$. We will always assume that $N$ is even.
The following construction yields the $osp(M|N)$ Yangian, and it
will lead to the \emph{non-super Yangians}  by taking $N = 0$, $\theta_0=1$
(orthogonal case) or $M=0$, $\theta_0=-1$ (symplectic case).
\begin{definition}
  For each index $i$, we introduce a sign $\theta_i$
  \begin{equation}
    \label{eq:deftheta}
    \theta_i = \begin{cases}
      +1 & \qmbox{for} 1\le i\le M+\frac{N}{2}  \cr
      -1 & \qmbox{for} M+\frac{N}{2}+1 \le i\le M+N
    \end{cases}
  \end{equation}
  and a conjugate index $\bar\imath$
  \begin{equation}
    \label{eq:defibar}
    \bar\imath =
    \begin{cases}
    M+1-i &\qmbox{for} 1\le i \le M \\
    2M+N+1-i &\qmbox{for} M+1\le i\le M+N
  \end{cases}
  \end{equation}
\end{definition}
In particular $\theta_{i} \theta_{\bar{\imath}} = \theta_{0} (-1)^{[i]} $.
\\
As usual $E_{ij}$ denotes the elementary matrix with entry 1 in row $i$
and column $j$ and zero elsewhere.
\begin{definition}
  For $A = \sum_{ij} \; A^{ij} \;E_{ij}$, we
  define the transposition $t$ by
  \begin{equation}
    \label{eq:t}
    A^t = \sum_{ij} (-1)^{[i][j]+[j]} \theta_i \theta_j \; A^{ij} \,
    E_{\bar\jmath \bar\imath}
    = \sum_{ij}  \left(A^{ij}\right)^t \,
    E_{ij}
  \end{equation}
  It satisfies $(A^t)^t = A$ and, for $\CC$-valued matrices,  $(AB)^t=B^tA^t$.
\end{definition}
We shall use a graded tensor product, i.e. such that, for $a$, $b$,
$c$ and $d$ with definite gradings,
$(a\otimes b)(c\otimes d) = (-1)^{[b][c]} ac \otimes bd$.
\begin{definition}
  Let $P$ be the (super)permutation operator
  (i.e.  $X_{21}\equiv PX_{12}P$)
  \begin{equation}
    \label{eq:Pdef}
    P = \sum_{i,j=1}^{M+N} (-1)^{[j]} E_{ij} \otimes E_{ji}
  \end{equation}
  and
  \begin{equation}
    \label{eq:Qdef}
    K \equiv
    P^{t_{1}} = \sum_{i,j=1}^{M+N} (-1)^{[i][j]} \theta_{i} \theta_{j}
    E_{\bar{\jmath}\bar{\imath}} \otimes E_{ji} \;,
  \end{equation}
  where $t_1$ is the transposition in the first space of the tensor
  product. In particular $P_{21}=P_{12}$ and $K_{21}=K_{12}$.
  \\
  We define the $R$-matrix
  \begin{equation}
    \label{eq:RPK}
    \displaystyle R(u) = \II + \frac{P}{u} - \frac{K}{u+\kappa}\;.
  \end{equation}
\end{definition}

\begin{proposition}
  The matrix $R(u)$ satisfies
  \begin{eqnarray}
    \label{eq:propR}
        && R_{12}^{t_{1}}(-u-\kappa) = R_{12}(u) \,, \qquad \mbox{(crossing
        symmetry)} \\
    \label{eq:Runit}
        && R_{12}(u) \, R_{12}(-u) = (1-1/u^2) \II \,, \qquad
        \mbox{(unitarity)}
  \end{eqnarray}
  provided that $2\kappa =
  (M-N-2)\theta_{0}=(\alpha_0+2\rho,\alpha_0)/2$, where $\rho$ is the
  super Weyl vector and $\alpha_0$ the longest root.
\end{proposition}
\textbf{Proof:} we use the fact that the operators $P$ and $K$ satisfy
  \begin{equation}
    \label{eq:PK}
    P^2 = \II, \qquad
    PK = KP = \theta_{0} K, \qquad \mbox{and} \qquad
    K^2 = \theta_0 (M-N)K
  \end{equation}

\finproof

\begin{theorem}
  The $R$-matrix (\ref{eq:RPK}) satisfies the super Yang--Baxter equation
  \begin{equation}
    \label{eq:YBE}
    R_{12}(u) \, R_{13}(u+v) \, R_{23}(v) = R_{23}(v) \, R_{13}(u+v)
    \, R_{12}(u)
  \end{equation}
  for $2\kappa = (M-N-2)\theta_{0}$, where the
  graded tensor product is understood.
\end{theorem}
\textbf{Proof:} we use the following relations obeyed by the matrices
$P$ and $K$
\begin{eqnarray}
  && P_{13} \, K_{23} = K_{12} \, P_{13} \qquad \qquad \qquad \quad
  K_{13} \, K_{12} = P_{23} \, K_{12} \nonumber\\
  && P_{12} \, P_{23} \, K_{12} = \theta_{0} \, P_{13} \, K_{12}
  \qquad \qquad P_{12} \, K_{23} \, K_{12} = \theta_{0} \,
  K_{13} \,
  K_{12} \nonumber\\
  && K_{12} \, K_{13} \, K_{23} = \theta_{0} \, P_{13} \, K_{23}
  \qquad \qquad 
  K_{12} \, P_{23} \, K_{12} = K_{12}
  \label{eq:relPK}
\end{eqnarray}
These relations are obtained by direct computation using the
definition of the matrices $P$ and $K$.

\finproof
\\
In the case related to $so(N)$,
this solution of the Yang--Baxter equation with spectral parameter was
found in \cite{ABZABZ78}. It is also one
of the cases explored in \cite{KS}. 

\section{Yangians\label{sect:yangians}}
\setcounter{equation}{0}

We consider the Hopf (super)algebra $\cU(R)$ generated by the operators
$T^{ij}_{(n)}$, for $1\le i,j\le M+N$, $n\in\ZZ_{\ge0}$, encapsulated
into a $(M+N)\!\times\!(M+N)$ matrix
\begin{equation}
  T(u)
  = \sum_{n \in \ZZ_{\ge 0}}
  T_{(n)} \, u^{-n} 
  = \sum_{i,j=1}^{M+N} \sum_{n \in \ZZ_{\ge 0}}
  T^{ij}_{(n)} \, u^{-n} \, E_{ij}
  = \sum_{i,j=1}^{M+N} T^{ij}(u) \, E_{ij}
  \label{eq:defT}
\end{equation}
and $T^{ij}_{(0)} = \delta_{ij}$.
One defines $\cU(R)$ by imposing the following constraints on
$T(u)$
\begin{equation}
  R_{12}(u-v) \, T_1(u) \, T_2(v) = T_2(v) \, T_1(u) \, R_{12}(u-v)
  \label{eq:RLL}
\end{equation}
with the matrix $R(u)$ defined in (\ref{eq:RPK}).
\\
The explicit commutation relations between the generating
operators $T^{ij}(u)$ read
\begin{eqnarray}
  \big[ T^{ij}(u) , T^{kl}(v) \big] &=&
  \frac{(-1)^{[k][i]+[k][j]+[i][j]}}{u-v}
  \; \Big(
  T^{kj}(v) T^{il}(u) - T^{kj}(u) T^{il}(v) \Big) \nonumber \\
  &+& \! 
  \; \frac{1}{u-v+\kappa} \ \sum_{p} \Big(
  \delta_{i\bar{k}} \;(-1)^{[p]+[j][i]+[j][p]}\;
  \theta_{\bar{\imath}}\theta_{\bar p} T^{pj}(u)
  T^{\bar{p}l}(v) \nonumber\\
  && \hspace*{72pt}
  -  \delta_{j\bar{l}}\; (-1)^{[k][j]+[i][k]+[i][p]}
  \;\theta_{\bar{p}}\theta_{\bar{\jmath}} T^{k\bar{p}}(v) T^{ip}(u) \Big)
  \label{eq:relcommTT}
\end{eqnarray}
that is, in terms of the generators $T^{ij}_{(n)}$
\begin{eqnarray}
&& \big[ T^{ij}_{(r+2)} , T^{kl}_{(s)} \big] + \big[ T^{ij}_{(r)} ,
T^{kl}_{(s+2)} 
\big] = 2 \big[ T^{ij}_{(r+1)} , T^{kl}_{(s+1)} \big] - \kappa \big[
T^{ij}_{(r+1)} , T^{kl}_{(s)} \big] + \kappa \big[ T^{ij}_{(r)} ,
T^{kl}_{(s+1)} 
\big] \nonumber \\
&& \hspace*{20mm} + \; (-1)^{[k][i]+[k][j]+[i][j]} \Big( T^{kj}_{(s)}
T^{il}_{(r+1)} - T^{kj}_{(r+1)} T^{il}_{(s)} 
- T^{kj}_{(s+1)} T^{il}_{(r)} 
\nonumber\\&& \hspace{60mm}
+ T^{kj}_{(r)} T^{il}_{(s+1)} + \kappa T^{kj}_{(s)} T^{il}_{(r)} -
\kappa T^{kj}_{(r)} 
T^{il}_{(s)} \Big) \nonumber \\
&& \hspace*{20mm} + \sum_{p} \Big( \delta_{i\bar{k}} \,
(-1)^{[p]+[j][i]+[j][p]} \, \theta_{\bar{\imath}}\theta_{\bar{p}} \,
(T^{pj}_{(r+1)} T^{\bar{p}l}_{(s)} - T^{pj}_{(r)}
T^{\bar{p}l}_{(s+1)}) \nonumber 
\\
&& \hspace*{35mm} - \delta_{j\bar{l}} \, (-1)^{[k][j]+[i][k]+[i][p]} \,
\theta_{\bar{p}}\theta_{\bar{\jmath}} \, (T^{k\bar{p}}_{(s)} T^{ip}_{(r+1)} -
T^{k\bar{p}}_{(s+1)} T^{ip}_{(r)}) \Big)
\label{eq:modes}
\end{eqnarray}
where $r,s\ge -2$ with by convention $T^{ij}_{(n)}=0$ for $n<0$.
\\
The Hopf algebra structure of $\cU(R)$ is given by \cite{FRT}
\begin{eqnarray}
  \label{eq:coproduct}
  && \Delta \big( T(u) \big) = T(u) \,\dot\otimes\, T(u) 
  \qquad \mbox{i.e.}\qquad
  \Delta \big( T^{ij}(u) \big) = \sum_{k=1}^{M+N} T^{ik}(u) \otimes
  T^{kj}(u) 
  \\
  && S(T(u)) =  T(u)^{-1}
  \quad ; \quad
  \epsilon(T(u)) = \II_{M+N}
  \label{eq:antipode}
\end{eqnarray}

\begin{theorem}
  The operators generated by $C(u) = T^t(u-\kappa) \, T(u)$ lie in the
  centre of the algebra $\cU(R)$ and $C(u) = c(u) \II$. 
  Furthermore,  $\Delta(c(u)) = c(u)\otimes c(u)$ and the two-sided
  ideal $\cI$ generated by 
  $C(u)-\II$ is also a coideal. The quotient $\cU/\cI$ is then a Hopf
  algebra. 
\end{theorem}
\textbf{Proof:} We first prove that $C(u)$ is diagonal.  Indeed, the
relation (\ref{eq:RLL}) implies
\begin{equation}
K_{12} \, T_1(u-\kappa) \, T_2(u) = T_2(u) \, T_1(u-\kappa) \, K_{12}
\label{eq:KTT}
\end{equation}
from which it follows, after having transposed in space 1,
\begin{equation}
\sum_{ijkl} (-1)^{[k]} T^t(u-\kappa)^{ij} T(u)^{jl} E_{ik}\otimes E_{kl} =
\sum_{pqsr} (-1)^{[p][s]+[p][r]+[s][r]} T(u)^{pq} T^t(u-\kappa)^{qr}
E_{sr}\otimes 
E_{ps}
\label{eq:TTEE}
\end{equation}
Therefore, one has
\begin{equation}
\sum_{j} T^t(u-\kappa)^{ij}  T(u)^{jl} = \delta_{il} \;c(u) \qquad
\mbox{or} \qquad C(u) = c(u) \, \II
\label{eq:TTf}
\end{equation}
Let us prove that $c(u)$ is a central element. One gets
\begin{eqnarray}
C(u) \, T_{2}(v) &=& T^{t}_1(u-\kappa) \, T_1(u) \, T_{2}(v) \nonumber \\
&=& T^t_1(u-\kappa) \, R_{12}^{-1}(u-v) \, T_{2}(v) \, T_1(u) \,
R_{12}(u-v)
\end{eqnarray}
where we have used the unitarity and crossing properties (\ref{eq:Runit})
and (\ref{eq:propR}) of $R(u)$.  Now using the transposition of the relation
(\ref{eq:RLL}) in space 1 and the crossing property of $R(u)$, one can derive
the following exchange relation:
\begin{equation}
T_1^{t}(u-\kappa) \, R_{12}^{-1}(u-v) \, T_2(v) = T_2(v) \,
R_{12}^{-1}(u-v) \, T_1^{t}(u-\kappa)
\label{eq:RLLtrans}
\end{equation}
Hence
\begin{eqnarray}
C(u) \, T_{2}(v) &=& T_{2}(v) \, R_{12}^{-1}(u-v) \, T^{t}_{1}(u-\kappa)
\, T_1(u) \, R_{12}(u-v) \nonumber \\
&=& T_{2}(v) \, R_{12}^{-1}(u-v) \, C(u) \, R_{12}(u-v)
\end{eqnarray}
Since $C(u) = c(u) \, \II$, one obtains easily $C(u) \, T_{2}(v) = T_{2}(v)
\, C(u)$.
\\
{}From the defining relations of $C(u)$ 
the coproduct of $c(u)$ is straightforwardly obtained as 
$
\Delta(c(u)) = c(u) \otimes c(u)
$
which shows that $\cI$ is a coideal. It is interesting to note that
this is precisely the structure of the coproduct of the quantum
determinant whenever such an object has been constructed.  
\finproof

\medskip
\noindent

At order $u^{-1}$ the equation $C(u)= \II$ yields the relation  
$T^t_{(1)} + T_{(1)} = 0$. Note that those linear relations 
${T^{ij}_{(1)}}^t + T^{ij}_{(1)} = 0$ for which $i\ne j$ were
already implied by the commutation relations (\ref{eq:modes}).
At higher orders,  $C(u)= \II$ induces relations with the generic form 
$T^t_{(n)} + T_{(n)} = \cF(T_{(m)},m<n)$ where $\cF$ is a quadratic
function. 
\\
In particular, once the exchange relations (\ref{eq:modes}) (for
$r=s=0$) are taken into account, the generators $T_{(1)}^{ij}$
exhibit the structure of the Lie (super) algebra $\fg$. 

\medskip
\noindent

\begin{definition}
  Let $\fg$ be a finite dimensional complex simple Lie (super)
  algebra.
  We define the bialgebra  
  $\fg[u]$ as $\fg \otimes_\CC \CC[u]$ 
  endowed with the Poisson cobracket $\delta$ defined by
  \begin{equation}
    \delta f(u,v) = 2 \left[\un\otimes f(v) + f(u)\otimes
    \un,\frac{\fC}{u-v}\right] \qquad 
    \label{eq:petitdelta}
  \end{equation}
  where $\fC$ is the tensorial Casimir element of $\fg$ associated
  with a given non-degenerate invariant bilinear form $\fB$,
  and $f : 
  \CC \to \fg$ is a 
  polynomial map, i.e. an element of $\fg[u]$.  
\end{definition}
\begin{theorem}
  Let $\fg$ be a finite dimensional complex simple Lie (super)
  algebra of type $so(M)$, $sp(N)$, $osp(M|N)$. 
  Let $\cU(R)$ be the Hopf algebra with 
  generators $T(u)$ subject to the relations (\ref{eq:RLL}) and Hopf
  structure (\ref{eq:coproduct})-(\ref{eq:antipode}).  The
  quotient of the algebra $\cU(R)$ by the two-sided ideal $\cI$ generated
  by $C(u) = T^t(u-\kappa) \, T(u) = \II$
  (i.e. $c(u) = 1$) is a homogeneous quantisation of $(\fg[u],\delta)$.
\end{theorem}
\textbf{Proof:} 
We define $\cU_\hbar$ as the algebra generated by the generating
functional 
$\tilde t(u)$ 
\begin{equation}
  \label{eq:ttilde}
  \tilde t(u) = \frac{1}{\hbar} \Big( T(u/\hbar) - 1 \Big)
\end{equation}
and the identity, the relations being derived from those of $\cU(R)$, i.e.
\begin{equation}
  \label{eq:ttildettilde}
  \begin{split}
    [\tilde t_1(u), \tilde t_2(v)] &= 
    \left[\tilde t_1(u) + \tilde t_2(v),\frac{P}{u-v}\right] 
    - \frac{\hbar}{u-v}(P\tilde t_1(u) \tilde t_2(v)
    -  \tilde t_1(u) \tilde t_2(v) P) \\
    &-
    \left[\tilde t_1(u) + \tilde t_2(v),\frac{K}{u-v+\hbar\kappa}\right] 
    +
    \frac{\hbar}{u-v+\hbar\kappa}(K\tilde t_1(u)  \tilde t_2(v) 
    -  \tilde t_1(u) \tilde t_2(v) K) \;.
  \end{split}
\end{equation}
Thus the relations in 
$\cU_\hbar/(\hbar\cU_\hbar)$ are 
\begin{eqnarray}
  \label{eq:demiboucle}
  [\tilde t_1(u), \tilde t_2(v)] &=& \left[\tilde t_1(u) + \tilde
  t_2(v),\frac{P-K}{u-v}\right] \;.
\end{eqnarray}
\\
The equation $C(u)= \II$ expressed in $\cU_\hbar$ generate a two-sided
ideal $\cI_\hbar$, which now
induces relations with the generic form 
$\tilde t^t_{(n)} + \tilde t_{(n)} = \hbar\, \cF(\tilde t_{(m)},m<n)$
where $\cF$ is a quadratic 
function. In the quotient algebra $\cU_\hbar/(\hbar\cU_\hbar)$ this
becomes equivalent to the standard linear symmetrisation relation 
$J^t_{(n)} + J_{(n)} =0$  for the generators of the loop algebra~$\fg[u]$,
so that $\cU'_\hbar/(\hbar\cU'_\hbar) \simeq \cU(\fg[u])$ as algebras,
for $\cU'_\hbar\equiv \cU_\hbar/{\cI_\hbar}$.  
This characterises $\cU'_\hbar$ as a quantisation of the algebra
$\cU(\fg[u])$. 
\\
We now examine the coproduct structure in order to recognise it as a
quantisation of the cocommutator $\delta$, namely 
\begin{equation}
  \label{eq:quantisation}
  \left. \frac{\Delta - \Delta^{op}}{\hbar}\; (\tilde
  t(u))\right|_{\!\!\mod \hbar} 
  = \delta\left(\tilde t(u) \big|_{\!\!\mod \hbar} \right)
\end{equation}
{}From (\ref{eq:coproduct}), the order $u^{-n}$ of the $(i,j)$ entry of the
left hand side of this formula reads
\begin{equation}
  \label{eq:DeltamoinsDeltaop}
  \left. \frac{\Delta - \Delta^{op}}{\hbar} \; \tilde t_{(m)}
  \right|_{\!\!\mod \hbar}   = 
  \sum_{r=0}^{m} 
  \left.\left(
  \tilde t_{(r)}\,\dot\otimes\, \tilde t_{(m-r)} 
  - \tilde t_{(r)} \,\dot\otimes\,   \tilde t_{(m-r)} 
  \right)\right|_{\!\!\mod \hbar}  \;.
\end{equation}
Now, denoting generically $\tilde t^a = \tilde t^{ij}-(\tilde t^{ij})^t$
and
$E_a = E_{ij}-(E_{ij})^t$,
and using $\tilde t = \tilde t^t \mod \hbar$, we can symmetrise and
get 
\begin{eqnarray}
  \sum_c 
  \left. \frac{\Delta - \Delta^{op}}{\hbar} \; \tilde t^c_{(m)} E_c
  \right|_{\!\!\mod \hbar}   &=& 
  \sum_{a,b}\sum_{r=0}^{m} 
  \tilde t^a_{(r)}\otimes \tilde t^b_{(m-r)} 
  \;[E_a,E_b]
  \;\bigg|_{\!\!\mod \hbar}  
  \nonumber\\
  &=&
  \sum_{a,b}\sum_{r=0}^{m} 
  \tilde t^a_{(r)}\otimes \tilde t^b_{(m-r)} 
  \; {f_{ab}}^c \, E_c \;
  \bigg|_{\!\!\mod \hbar}  
\end{eqnarray}
The right hand side of the formula (\ref{eq:quantisation})
can be computed once one recalls
that $\fC = \sum_{ab} \fB_{ab} t^a \otimes t^b$.
One obtains 
\begin{equation}
  \delta(t^a_{(m)}) =   \sum_{a,b}\sum_{r=0}^{m} 
   t^a_{(r)}\otimes t^b_{(m-r)} 
  \; {f^c}_{ab} 
\end{equation}
where the $t^a_{(m)}$ denote the generators of the loop algebra
$\fg[u]$.
Since the modes of $(\tilde t(u) \big|_{\!\!\mod \hbar})$ coincide
with the $t^a_{(m)}$ and the structure constants ${f_{ab}}^c$ and
${f^c}_{ab}$ are identified through the bilinear form $\fB$, one gets 
the desired result (\ref{eq:quantisation}).
\\
Therefore the Hopf algebra $\cU(R)/\cI \equiv \cU'_{\hbar=1}$ is a
quantisation of $\cU(\fg[u])$ and $\Delta$ is a quantisation of
$\delta$. \\
\null\hfill\finproof

\null
\medskip\noindent
{}From the above theorem, we are naturally led to the following
definition: 
\begin{definition}
We define the Yangian of $osp(M|N)$ as   
$Y(\fg)  \equiv \cU(R)/\cI$.\\
Explicitly, its defining relations are given by
\begin{eqnarray*}
&& R_{12}(u-v)\; T_1(u)\; T_2(v) = T_2(v)\; T_1(u)\; R_{12}(u-v) \;,\\
&& C(u) = T^t(u-\kappa) \, T(u) = \II\,,
\end{eqnarray*}
where 
$\displaystyle R_{12}(u)=\II + \frac{P}{u} - \frac{K}{u+\kappa}$\;.
\end{definition}
For $N=0$ or $M=0$,  this
definition is consistent with the one of Drinfel'd \cite{Dri85,Dri86} 
for the $so(M)$ and $sp(N)$ cases respectively.

\medskip
\noindent
\textbf{Remark:} The explicit $R$-matrices for the
Yangians $Y(so(N))$ and 
$Y(sp(N))$ can be obtained by taking the scaling limit $q \to 1$, $z = q^u
\to 1$ keeping $u$ fixed, of the evaluated trigonometric $R$-matrices of
$\cU_{q}(\widehat{so(N)})$ and $\cU_{q}(\widehat{sp(N)})$ computed in
\cite{Jim86}. 
Similarly, one can show that the $R$-matrix of $Y(osp(1|2))$ is the
scaling limit of the evaluated trigonometric $R$-matrix of
$\cU_{q}(\widehat{osp(1|2)})$ \cite{KTtwist}.

\section{Twisted Yangians and reflection algebras}
\setcounter{equation}{0}
We would finally like to comment upon a possible connection
between the notions of twisted Yangians and reflection algebras
within the framework of this Yangian construction. 
\\
Following the lines of \cite{olsh,MNO} (see also \cite{briot} 
for the supersymmetric case), we define on $\cU(R)$:
\begin{equation}
\tau[T(u)]=T^t(-u-\kappa)
\end{equation}
which reads for the super-Yangian generators:
\begin{equation}
\tau(T^{ab}(u))=(-1)^{[a]([b]+1)}\theta_{a}\theta_{b}\ 
T^{\bar b\bar a}(-u-\kappa)
\label{supertau}
\end{equation}
$\tau$ is  an algebra automorphism, as a
direct consequence of unitarity, crossing symmetry and the property 
$R^{t_1t_2}(u)=R(u)$ which itself comes from $P^{t_1t_2}=P$.
\\
The twisted super-Yangian $\cU(R)^{tw}$ is the subalgebra generated 
by $S(u)=\tau[T(u)]T(u)$, with $\tau$ given in (\ref{supertau}).
$S(u)$ obeys the following relation:
\begin{equation}
  R_{12}(u-v)\, S_{1}(u)\, R_{12}(u+v)\, S_{2}(v) = 
  S_{2}(v)\, R_{12}(u+v)\, S_{1}(u)\, R_{12}(u-v)
  \label{rsrs}
\end{equation}
It is easy to show that $\cU(R)^{tw}$ is a coideal in $\cU(R)$.
\\
Similarly, one introduces the notion of reflection algebras
$\cS(R)$, generated by 
\begin{equation}
  B(u)=T^{-1}(-u)T(u)
\end{equation}
which obeys the same relation (\ref{rsrs}), interpreted here as a
reflection equation. $\cS(R)$ is also a coideal of $\cU(R)$.
This type of algebras have been originally introduced in \cite{skly}
for the Yangian $Y(N)$, based on $gl(N)$, and play an important role
in integrable systems with boundaries (see e.g. \cite{BNLS}).
\\
However in the coset $\cU(R)/\cI$, one has $B(u)=S(u)$, so that
$\cS(R)$ and $\cU(R)^{tw}$ are two versions of 
the same Hopf coideal in $\cU(R)$. The situation is here different 
from the case of the Yangian $Y(N)$.
Indeed the twisted Yangians $Y^\pm(N)$ and 
the boundary algebras $B(N,\ell)$ are known to be different for $N>2$, 
whilst for $N=2$ one has $B(2,0)=Y^-(2)$ and $B(2,1)=Y^+(2)$
\cite{MR}.

\bigskip\noindent
\textbf{Acknowledgements:} We would like to thank A.~Molev and
V.~Tolstoy for discussions and comments.


\begin{thebibliography}{99}

\bibitem{Dri85}
  V.G. Drinfel'd, \textsl{Hopf algebras and the quantum Yang--Baxter
    equation,} 
  Soviet.  Math.  Dokl.  \textbf{32} (1985) 254--258.
  
\bibitem{Dri86}
  V.G. Drinfel'd, \textsl{Quantum Groups,}
  Proceedings Int. Cong. Math. Berkeley, California, USA (1986) 798--820.
  
\bibitem{Dri88}
  V.G. Drinfel'd, \textsl{A new realization of {Y}angians and quantized affine
    algebras,} Soviet.  Math.  Dokl.  \textbf{36} (1988) 212--216.
  
\bibitem{Molev2001}
  A.I. Molev, \emph{Yangians and their applications}, Handbook of
  {A}lgebra, vol. 3, Elsevier, to appear.

\bibitem{FRT} L.D.~Faddeev, N.Yu.~Reshetikhin and L.A.~Takhtajan,
  \textsl{Quantization of Lie groups and Lie algebras,} Leningrad Math.  J.
  \textbf{1} (1990) 193--225.
  
\bibitem{Naz91} M.L.~Nazarov, \textit{Quantum Berezinian and the
    classical Capelli identity,} Lett. Math. Phys. \textbf{21} (1991)
    123--131. 

\bibitem{Zha95} R.B.~Zhang, \textit{The $gl(M|N)$ super Yangian and
    its finite dimensional representations,}
    Lett. Math. Phys. \textbf{37} (1996) 419--434.

\bibitem{olsh}G.I. Olshanski, {\it Twisted Yangians and infinite dimensional 
    Lie algebras}, in ``Quantum groups'', Lecture Notes in Math. {\bf 1510} 
  (P. Kulish ed.), pp. 104--120, NY 1992.
  
\bibitem{MNO} A. Molev, M. Nazarov and G. Olshanski,
  \textsl{Yangians and classical {L}ie algebras,}
  Russ. Math. Surveys \textbf{51} (1996) 205--282,
  \texttt{hep-th/9409025}. 
  
\bibitem{briot} C. Briot and E. Ragoucy, 
  {\it Twisted superYangians 
    and their representations}, preprint LAPTH-875/01,
  \texttt{math.QA/0111308}.

\bibitem{ABZABZ78} Al.~B.~Zamolodchikov, Al.~B.~Zamolodchikov,
  \textsl{Relativistic factorized $S$-matrix in two dimensions having $O(N)$
  isotropic symmetry,} Nucl.  Phys.  \textbf{B133} (1978) 525--535 and
  \textsl{Factorized $S$-matrices in two dimensions as the exact solutions of
  certain relativistic quantum field models,} Ann.  Phys.  \textbf{120}
  (1979) 253--291.

\bibitem{KS} P.P.~Kulish, E.K.~Sklyanin,
  \textsl{Solutions of the Yang--Baxter equation,}
  Zap. Nauchn. Sem. LOMI, \textbf{95} (1980) 129--160 and
  J. Sov. Math. \textbf{19} (1982) 1596--1620.
  
\bibitem{Isa} A.P. Isaev,
  \textsl{Quantum groups and Yang--Baxter equations,} Phys. Part. Nucl.
  \textbf{26} (1995) 501--526.
  
\bibitem{KR} A.N. Kirilov and N.Yu. Reshetikhin,
  \textsl{The Yangians, Bethe Ansatz and combinatorics,}
  Lett. Math. Phys. \textbf{12} (1986) 199-208.
  
\bibitem{Jim86}
  M.~Jimbo, \textsl{Quantum {$R$}-matrix for the generalized {T}oda
    system}, Commun.  Math.  Phys.  \textbf{102} (1986), 537--547.

\bibitem{KTtwist}
  S.~Khoroshkin and V.~Tolstoy, {\it Twisting of quantum (super)algebras.
    {C}onnection of {D}rinfeld's and {C}artan--{W}eyl realizations for quantum
    affine algebras}, \texttt{hep-th/9404036}.
  
\bibitem{skly} E.K. Sklyanin, {\it Boundary conditions for integrable 
    quantum systems}, J. Phys. {\bf A21} (1988) 2375--2389.
  
\bibitem{BNLS} M. Mintchev, E. Ragoucy and P. Sorba, {\it Spontaneous 
    symmetry breaking in the gl(N)-NLS algebra}, J. Phys. 
  {\bf A34} (2001) 8345--8364, 
  {\tt hep-th/0104079}.
  
\bibitem{MR} A. Molev and E. Ragoucy, {\it Representations of boundary
    algebras}, Rev. Math. Phys. {\bf 14} (2002)  317--342, 
{\tt math.QA/0107213}.
  
\end{thebibliography}
\end{document}